\documentclass[letterpaper, 10 pt, conference, twoside]{IEEEtran}
\usepackage{fancyhdr}
\usepackage{graphicx,subfigure}         
\usepackage{url}
\usepackage{amsmath}
\usepackage{amssymb}
\usepackage{amsthm} 
\usepackage{lipsum,adjustbox}
\usepackage{multirow}
\let\labelindent\relax
\usepackage[shortlabels]{enumitem}
\usepackage{mathtools}
\usepackage{cite}
\usepackage{color}
\allowdisplaybreaks
\usepackage[mathscr]{eucal}
\DeclareMathOperator{\D}{\mathcal{D}}
\DeclareMathOperator{\V}{\mathcal{V}}
\DeclareMathOperator{\N}{\mathcal{N}}
\DeclareMathOperator{\NN}{\mathscr{N}}
\DeclareMathOperator{\E}{\mathcal{E}}
\DeclareMathOperator{\DD}{\mathcal{D}}
\DeclareMathOperator{\G}{\mathcal{G}}
\newcommand{\vv}{\mathbf{v}}
\newcommand{\zz}{\mathbf{z}}
\newcommand{\PP}{\mathbf{P}}
\newcommand{\qq}{\mathbf{q}}
\newcommand{\bb}{\mathbf{b}}
\newcommand{\ww}{\mathbf{w}}
\newcommand{\yy}{\mathbf{y}}
\newcommand{\real}{\mathbb{R}}
\newcommand{\naturals}{\mathbb{N}}
\newcommand{\diag}{\operatorname{diag}}
\newcommand{\one}{\mathbf{1}}
\newcommand{\zero}{\mathbf{0}}
\newcommand{\nulls}{\operatorname{null}}
\newcommand{\grad}{\psi_{\text{grad}}}
\newcommand{\est}{\psi_{\text{est}}}
\newcommand{\PPest}{\PP_{\text{est}}}
\newcommand{\ic}{\psi_{\text{interc}}}
\newcommand{\linear}{\phi_{\text{grad}}}
\newcommand{\revision}[1]{{#1}}

\newtheorem{theorem}{Theorem}[section]
\newtheorem{proposition}[theorem]{Proposition}

\theoremstyle{remark}
\newtheorem{remark}{Remark}
\theoremstyle{definition}

\newcommand{\longthmtitle}[1]{\mbox{} \textit{(#1):}}
\newcommand{\until}[1]{\{1,\dots,#1\}}
\newcommand{\setdef}[2]{\{#1 \; | \; #2\}}

\newcommand{\map}[3]{#1:#2 \rightarrow #3}

\newcommand{\oprocendsymbol}{\hbox{$\bullet$}}
\newcommand{\oprocend}{\relax\ifmmode\else\unskip\hfill\fi\oprocendsymbol}
                               
\graphicspath{ {Figures/} }
  
\title{Network Optimization via Smooth Exact Penalty
  Functions Enabled by Distributed Gradient Computation \thanks{This
    work was supported by the ARPA-e NODES program, Cooperative
    Agreement DE-AR0000695 and NSF Award ECCS-1917177. A preliminary
    version of this paper appeared at the IEEE Conference on Decision
    and Control as~\cite{PS-JC:18-cdc}.}
}

\author{Priyank Srivastava \quad Jorge
  Cort\'{e}s
  \thanks{The authors are with the Department of Mechanical and
    Aerospace Engineering, UC San Diego,
    {\tt\small \{psrivast,cortes\}@ucsd.edu}}%
}

\graphicspath{ {Figures/} }
  
\begin{document}
  
 \maketitle
\thispagestyle{fancy}
  
  \begin{abstract}
    This paper proposes a distributed algorithm for a network of
    agents to solve an optimization problem with separable objective
    function and locally coupled constraints. Our strategy is based on
    reformulating the original constrained problem as the
    unconstrained optimization of a smooth (continuously
    differentiable) exact penalty function.  Computing the gradient of
    this penalty function in a distributed way is challenging even
    under the separability assumptions on the original optimization
    problem.  Our technical approach shows that the distributed
    computation problem for the gradient can be formulated as a system
    of linear algebraic equations defined by separable problem data.
    To solve it, we design an exponentially fast, input-to-state
    stable distributed algorithm that does not require the individual
    agent matrices to be invertible.  We employ this strategy to
    compute the gradient of the penalty function at the current
    network state. Our distributed algorithmic solver for the original
    constrained optimization problem interconnects this estimation
    with the prescription of having the agents follow the resulting
    direction. Numerical simulations illustrate the convergence and
    robustness properties of the proposed algorithm.
  \end{abstract}
  

\section{Introduction}
Network optimization problems arise naturally as a way of encoding the
coordination task entrusted to a multi-agent system in many areas of
engineering, including power, communication, transportation, and swarm
robotics.  The large-scale nature of these network problems together
with technological advances in communication, embedded computing, and
parallel processing have sparked the development of distributed
algorithmic solutions that scale with the number of agents, provide
plug-and-play capabilities, and are resilient against single points of
failure.  This paper is a contribution to the growing body of work
that deals with the design and analysis of provably correct
distributed algorithms that solve constrained optimization
problems with separable objective functions and locally expressible
constraints.  The novelty of our approach lies in the use of
continuously differentiable exact penalty functions to deal with the
constraints, thereby avoiding the characteristic chattering behavior
associated with non-differentiable approaches, and the reliance on
gradient descent directions, thereby avoiding the oscillatory behavior
characteristic of primal-dual schemes.

\textit{Literature review:} The breadth of applications of distributed
convex optimization~\cite{DPB-JNT:97,MGR-RDN:05,PW-MDL:09} has
motivated a growing body of work that builds on consensus-based
approaches to produce rich algorithmic designs with asymptotic
convergence guarantees, see~\cite{AN:14-sv} for a comprehensive
survey.  In this class of problems, each agent in the network
maintains, communicates, and updates an estimate of the complete
solution vector, whose dimension is independent of the network size.
This is in contrast to the setting considered here, where the
structure of the optimization problem lends itself to having instead
each agent optimize over and communicate its own local variable.
Considered collectively, these variables give rise to the solution
vector. Distributed algorithms to address this setting fall under
Lagrangian-based approaches that rely on primal-dual updates,
e.g.,~\cite{DF-FP:10,DR-JC:15-tac,EM-CZ-SL:17,SB-NP-EC-BP-JE:11,YX-TH-KC-ZL-GY-MF:17,SAA-KY-AHS:19} or
unconstrained reformulations that employ non-smooth penalty
functions~\cite{LX-SB:06,BJ-MJ:09,AC-JC:15-tcns}.
Our approach here is based on the exact reformulation of the original
problem using continuously differentiable penalty
functions~\cite{TG-EP:79,GdP-LG:89,SL:92,GDP:94}.
The work~\cite{GdP-LG:89} establishes, under appropriate regularity
conditions on the feasibility set, the complete equivalence between
the solutions of the original constrained and the reformulated
unconstrained optimization problems. The work~\cite{SL:92} proposes a
continuously differentiable exact penalty function that relaxes some
of the assumptions of~\cite{GdP-LG:89}. Notably, the works on
continuously differentiable exact penalty functions use centralized
optimization algorithms because the computations involved in the
definition of the unconstrained penalty function are of a centralized
nature. Our recent work~\cite{PS-JC:21-csl} provides a framework to
extend Nesterov acceleration to constrained optimization by
investigating conditions under which the penalty function is convex.

\textit{Statement of contributions:} We consider nonlinear programming
problems with a separable objective function and locally coupled
constraints. The starting point for our algorithm design is the exact
reformulation of the problem as an unconstrained optimization of a
continuously differentiable exact penalty function. Motivated by
enabling the computation of the gradient of this function by the
network agents, our first contribution is the design of a distributed
algorithm to solve a system of linear algebraic equation whose
coefficient matrix and constant vector can be decomposed as the
aggregate of (not necessarily invertible) coefficient matrices and
constant vectors, one per agent.
We establish the exponential convergence and characterize the
input-to-state stability properties of this algorithm. Building on it,
our second contribution is the structured computation of the gradient
of the penalty function in a distributed way. We accomplish this by
showing that the calculation of certain non-distributed terms in the
gradient can be formulated as solving appropriately defined systems of
linear algebraic equations defined by separable data.  Our third and
last contribution is the design of the distributed algorithm that
solves the original constrained optimization problem. This algorithm
is based on following gradient descent of the penalty function while
estimating the actual value of the gradient with the distributed
strategy that solves systems of linear algebraic equations.  We
establish the convergence of the resulting interconnection and
illustrate its performance in simulation, comparing it with
alternative approaches. We end by noting that, since the proposed
approach relies on the distributed computation of the gradient, the
methodology can also be used for accelerated distributed optimization
using Nesterov's method, something which we also illustrate 
in simulation.

\section{Preliminaries}\label{sec:prelim}
In this section, we present our notational conventions and review
basic concepts on graph theory and constrained optimization.

\textit{Notation:} Let $\real$ and $\naturals$ be the set of real and
natural numbers, resp. We let $\mathcal{X}^o$ and
$\overline{\mathcal{X}}$ denote the interior and closure
of~$\mathcal{X}$, resp.  For a real-valued function $f: \mathbb{R}^n
\rightarrow \mathbb{R}$, we let $\nabla f$ denote its gradient. When
we take the partial derivative with respect to a specific
argument~$x$, we employ the notation $\nabla_x f$. We denote vectors
and matrices by lowercase and uppercase letters, respectively. With a
slight abuse of notation, we let $(a;b)$ denote the concatenated
vector containing the entries of vectors $a$ and $b$, in that order.
$A'$ denotes the transpose of a matrix $A$. $A\otimes B$ denotes the
Kronecker product of two matrices $A$ and $B$. We use $\mathbf{0}$ and
$\one$ to denote the vector or matrix of zeros and ones of appropriate
dimension, respectively.  $\diag (v) \in \real^{n \times n}$ denotes
the diagonal matrix with the elements of $v \in \real^n$ in its
diagonal. Similarly, for a group of square matrices $\{A_i\}_{i \in
  \until{n}} \in \real^{m \times m}$, $\diag (A_i) \in \real^{mn
  \times mn}$ denotes the block-diagonal matrix with each of the
matrices $A_i$ arranged along the principal diagonal. We use
$\lambda_2(A)$ to denote the smallest non-zero eigenvalue of matrix
$A$, regardless of the multiplicity of eigenvalue $0$. null$(A)$
denotes the nullspace or kernel of a matrix $A$. We use dim$(W)$ to
denote the dimension of vector space $W$.

\textit{Graph theory:} We present basic concepts from graph theory
following~\cite{CDG-GFR:01}.  We denote an undirected graph by
$\G=(\V,\E)$, with $\mathcal{V}$ as the set of vertices and
$\E\subseteq\mathcal{V} \times \V$ as the set of edges. $(i,j) \in \E$
if and only if $(j,i) \in \E$. A vertex $j \in \V$ is a neighbor of
$i$ iff $(i,j) \in\E$, and $k$ is a 2-hop neighbor of $i$ if there
exists $j \in \mathcal{V}$ such that $(i,j) \in \E$ and
$(j,k) \in \E$. 
The set of all 1-hop neighbors of $i$ is denoted by $\mathscr{N}_i$.
A graph is  connected if there exists a path
between any two vertices. The degree of a node is the number of edges
connected to it. The degree matrix $D \in \real^{n \times n}$ is the
diagonal matrix with $D_{ii}=\text{deg}(v_i)$. The adjacency matrix
$A \in \real^{n \times n}$ is defined by $A_{ij}=1$ if
$(i,j) \subseteq \E$ and $A_{ij}=0$ otherwise. The Laplacian matrix is
$L=D-A$. Note that $\mathbf{1}'L=0$ and $0$ is a simple eigenvalue of
$L$ if and only if the graph $\G$ is connected.

\textit{Constrained optimization:} Here, we introduce basic concepts
of constrained optimization following~\cite{DPB:99}.  Consider the
following nonlinear optimization problem
\begin{equation}\label{eq:nl}
  \begin{aligned}
    & \min_{x \in \mathcal{D}}
    & & f(x) \\
    &\; \; \text{s.t.}  & & g(x) \leq 0 , \; h(x)=0,
  \end{aligned}
\end{equation}
where
$f:\real^n \rightarrow \real, \; g:\real^n \rightarrow \real^m, \; h:
\real^n \rightarrow \real^p$ are twice continuously differentiable
functions with $p \leq n$ and $\mathcal{D} \subset \real^n$ is a
compact set which is regular (i.e., $\DD= \overline{\mathcal{D}^o}$).
The feasible set of~\eqref{eq:nl} is
$\mathcal{F}= \setdef{x}{x \in \mathcal{D}, g(x) \leq0, h(x)=0}$.
Based on the index sets for the inequality constraints
\begin{align*}
  I_0(x)=\setdef{j}{g_j(x)=0} ,
  \\
  I_+(x)=\setdef{j}{g_j(x) \geq 0} ,
\end{align*}
we define the following regularity conditions:
\begin{enumerate}[(a)]
\item The \emph{linear independence constraint qualification} (LICQ)
  holds at $x \in \real^n$ if
  $\{\nabla g_j(x) \}_{j \in I_0(x)} \cup \{\nabla h_k\}_{k \in
    \until{p}}$ are linearly independent;
\item The \emph{extended Mangasarian-Fromovitz constraint
    qualification} (EMFCQ) holds at $x \in \real^n$ if
  $\{ \nabla h_k \}_{k \in \until{p}}$ are linearly independent and
  there exists $z\in \real^n$ with
  \begin{subequations}
    \begin{align}\label{mfcq}
      \nabla g_j(x)'z &< 0, \quad \forall j \in I_+(x) ,
      \\
      \nabla h_k(x)'z &=0, \quad \forall k \in \until{p} .
    \end{align}
  \end{subequations}
\end{enumerate}

The Lagrangian function
$L:\real^n \times \real^m \times \real^p \rightarrow \real$ associated
with \eqref{eq:nl} is given by
\begin{align*}
  L(x,\lambda,\mu) = f(x)+\lambda'g(x)+\mu'h(x) ,
\end{align*}
where $\lambda \in \real^m$ and $\mu \in \real^p$ are the Lagrange
multipliers (also called dual variables) associated with the
inequality and equality constraints, resp.  A Karush-Kuhn-Tucker (KKT)
point for~\eqref{eq:nl} is a triplet
$(\bar{x},\bar{\lambda},\bar{\mu})$ such that
\begin{align*}
  \nabla_xL(\bar{x},\bar{\lambda},\bar{\mu})
  &=0 ,
  \\
  \bar{\lambda}'g(\bar{x}) 
  &=0 , \quad \bar{\lambda} \geq 0 , \quad
    g(\bar{x}) \leq0 ,
  \\
  \quad h(\bar{x})&=0 .
\end{align*}
Under any of the regularity conditions above, the KKT conditions are
necessary for a point to be locally optimal.

\textit{Continuously differentiable exact penalty functions:} With
exact penalty functions, the basic idea is to replace the
constrained optimization problem~\eqref{eq:nl} by an equivalent
unconstrained problem.  Here, we introduce continuously differentiable
exact penalty functions following~\cite{TG-EP:79,GdP-LG:89}.
Beyond the knowledge of the availability of such functions,
  the reader can defer parsing through the specific technical details
  below until they become critical in
  Section~\ref{sec:dist_formulation} below.  The key observation is
that one can interpret a KKT tuple as establishing a relationship
between a primal solution $\bar{x}$ and the dual variables
$(\bar{\lambda},\bar{\mu})$.  In turn, the following result introduces
multiplier functions that extend this relationship to any $x\in
\real^n$.

\begin{proposition}\longthmtitle{Multiplier functions and their
    derivatives~\cite{GdP-LG:89}}\label{prop:lambda}
  Assume that LICQ is satisfied at all $x \in \mathcal{D}$. Let
  $G(x)=\diag(g(x))$ and, for $\gamma \neq 0$, define $N:\real^n
  \rightarrow \real^{(m+p) \times (m+p)}$ by
  \begin{align}\label{eq:n}
    N(x)=
    \begin{bmatrix}
      \nabla g(x)'\nabla g(x) + \gamma^2 G^2(x) & \nabla g(x)' \nabla
      h(x)
      \\
      \nabla h(x)' \nabla g(x) & \nabla h(x)' \nabla h(x)
    \end{bmatrix} .
  \end{align}
  Then $N(x)$ is a positive definite matrix for any $x \in
  \D$. Given the functions $x \mapsto (\lambda(x), \mu(x))$
  defined by
  \begin{align}\label{eq:lambda}
    \begin{bmatrix}
      \lambda(x) \\ \mu(x)
    \end{bmatrix}= -N^{-1}(x) \begin{bmatrix} \nabla g(x)' \\ \nabla
      h(x)'
    \end{bmatrix} \nabla f(x) ,
  \end{align}
  one has that
  \begin{enumerate}[(a)]
  \item if $(\bar{x},\bar{\lambda},\bar{\mu})$ is a KKT triple for
    problem~\eqref{eq:nl}, then $\lambda(\bar{x})=\bar{\lambda}$ and
    $\mu(\bar{x})=\bar{\mu}$;
  \item both functions are continuously differentiable and their
    Jacobian matrices are given by
    \begin{align}\label{eq:grad_lambda}
      \begin{bmatrix}
        \nabla \lambda(x)' \\ \nabla\mu(x)'
      \end{bmatrix}= - N^{-1}(x) \begin{bmatrix} R(x) \\ S(x)
      \end{bmatrix} ,
    \end{align}
    where
    \begin{subequations}\label{eq:rs}
      \begin{align}\label{eq:r}
        R(x)
        &= \nabla g(x)' \nabla^2_x L(x,\lambda(x),\mu(x)) \notag
        \\
        & \quad +\sum\limits_{j=1}^m e_j^m\nabla_x L(x,\lambda(x),\mu(x))'
          \nabla^2 g_j(x)
        \\
        & \quad + 2\gamma^2 \Lambda(x) G(x) \nabla g(x)'
          \notag
        \\
        \label{eq:s}
        S(x) &= \nabla h(x)' \nabla^2_x
               L(x,\lambda(x),\mu(x)) \notag
        \\
        & \quad +\sum\limits_{k=1}^p
          e_k^p\nabla_x L(x,\lambda(x),\mu(x))' \nabla^2 h_k(x)
      \end{align}
    \end{subequations}
    where we use the shorthand notation
    \begin{align*}
      \nabla_x L(x,\lambda(x),\mu(x))
      &= [\nabla_x
        L(x,\lambda,\mu)]_{\substack{\lambda=\lambda(x) \\ \mu=\mu(x)}}
      ,
      \\
       \nabla^2_x L(x,\lambda(x),\mu(x)) &= [\nabla^2_x
                                           L(x,\lambda,\mu)]_{\substack{\lambda=\lambda(x) \\ \mu=\mu(x)}}
      ,
    \end{align*}
    $ \Lambda(x) =\diag(\lambda(x))$, and $e_j^m$ and $e_k^p$ denote,
    resp., the $j$th and $k$th column of the $m \times m$ and
    $p \times p$ identity matrix.
  \end{enumerate}
\end{proposition}

The multiplier functions in Proposition~\ref{prop:lambda} can be used
to replace the multiplier vectors in the augmented Lagrangian
of~\cite{RTR:74} to define the continuously differentiable exact
penalty function. Given $\epsilon > 0$ and $j \in \until{m}$, define
\begin{align*}
  y_j^{\epsilon}(x)
  &= \Big(-\min
    \Big[0,g_j(x)+\frac{\epsilon}{2}\lambda_j(x) \Big] \Big)^{1/2} ,
\end{align*}
and let $Y^{\epsilon}(x)=\diag(y^{\epsilon}(x))$.  Consider the
continuously differentiable function
$\map{f^\epsilon}{\real^n}{\real}$,
\begin{align}\label{eq:penalty}
  f^{\epsilon}(x)
  &=f(x)+\lambda(x)'(g(x)+Y^{\epsilon}(x)y^{\epsilon}(x))+\mu(x)'h(x)\notag
  \\
  &\quad +\frac{1}{\epsilon}\|g(x)+Y^{\epsilon}(x)y^{\epsilon}(x)\|^2 +
  \frac{1}{\epsilon}\|h(x)\|^2 .
\end{align}
The following result characterizes the extent to which $f^\epsilon$ is
an exact penalty function.

\begin{proposition}\longthmtitle{Continuously differentiable exact
    penalty function~\cite{GdP-LG:89}}\label{prop:exactness}
  Assume LICQ is satisfied at all $x \in \mathcal{D}$ and consider
  the unconstrained problem
  \begin{align}\label{eq:unc}
    \min_{x \in \mathcal{D}^o} f^\epsilon(x) .
  \end{align}
  Then, the following holds:
  \begin{enumerate}[(a)]
  \item there exists $\bar{\epsilon}$ such that the set of global
    minimizers of~\eqref{eq:nl} and~\eqref{eq:unc} are equal for all
    $\epsilon \in (0,\bar{\epsilon}]$;
  \item if $(\bar{x},\bar{\lambda},\bar{\mu})$ is a KKT point for
    problem~\eqref{eq:nl}, then $\nabla f^\epsilon(\bar{x})=0$ for all
    $\epsilon > 0$;
  \item under the additional assumption that EMFCQ holds on
    $\mathcal{D}$, there exists $\bar{\epsilon}$ such that for all
    $\epsilon \in (0,\bar{\epsilon}]$, $\nabla f^\epsilon(\bar{x})=0$
    implies that $(\bar{x},\lambda(\bar{x}),\mu(\bar{x}))$ is a KKT
    point for problem~\eqref{eq:nl}.   \end{enumerate}
\end{proposition}

Given the result of Proposition~\ref{prop:exactness}, we next turn our
attention to solve the unconstrained optimization
problem~\eqref{eq:unc}.
The next result, whose proof is given in the appendix, characterizes
the extent to which the gradient descent dynamics of $f^\epsilon$
satisfies the constraints while finding the optimizers of the original
constrained optimization problem.

\begin{proposition}\longthmtitle{Constraint satisfaction under
    gradient dynamics of penalty function}\label{prop:feasible}
  Given the optimization problem~\eqref{eq:nl}, assume LICQ is
  satisfied at all $x \in \mathcal{D}$. Consider the gradient dynamics
  $ \dot{x}=-\nabla f^\epsilon (x)$ of the penalty function
  $f^\epsilon$ in~\eqref{eq:penalty}.  Then, if at any time $t_0$,
  $x(t_0) \! \in \! \mathcal{F}$, we have
  \begin{enumerate}
  \item \textit{(Equality constraints):} $ x(t) \in \mathcal{F}$, for
    all $ t \geq t_0$ and all $\epsilon>0$ if the
    problem~\eqref{eq:nl} has just equality constraints;
  \item \textit{(Scalar inequality constraint):} there exists
    $\bar{\epsilon}>0$ such that $x(t) \in \mathcal{F}$, for all $t
    \geq t_0$ and all $\epsilon \in (0, \bar{\epsilon}]$ if the
    problem~\eqref{eq:nl} has only one inequality constraint;
  \item \textit{(General constraints):} in general, there is no guarantee that
    the evolution of the gradient dynamics stays feasible when the
    problem~\eqref{eq:nl} has more than one constraint if one of them
    is an inequality.
  \end{enumerate}
\end{proposition}

\section{Problem Statement}\label{sec:problem}
We consider separable network optimization problems where the overall
objective function is the aggregate of individual objective functions,
one per agent, and the constraints are locally expressible. Formally,
consider a group of $n \in \naturals$ agents whose interaction is
modeled by an undirected connected graph $\mathcal{G}=(\V,\E)$.  Each
agent $i \in \V$ is responsible for a decision variable
$x_i \in \real$.  Agent $i$ is equipped with a twice continuously
differentiable function $f_i:\real \rightarrow \real$.  The
optimization problem takes the form
\begin{equation}\label{eq:nlp}
  \begin{aligned}
    & \min_{x \in \mathcal{D}}
    & & f(x)=\sum\limits_{i=1}^n f_i(x_i) \\
    &\; \; \text{s.t.}  & & g(x) \leq 0 , \; h(x)=0,
  \end{aligned}
\end{equation}
with twice continuously differentiable vector-valued functions
$g:\real^n \rightarrow \real^m$, $h:\real^n \rightarrow \real^p$, and
$p \leq n$. Each component $\{\map{g_j}{\real^n}{\real}
  \}_{j=1}^{m}$ and $\{\map{h_k}{\real^n}{\real} \}_{k=1}^{p}$ of the
constraint functions is locally expressible.  Such kind of
  coupled constraints arise in numerous applications, such as
  power~\cite{ED-HZ-GBG:13}, communication~\cite{FPK-AKM-DKHT:98}, and
  transportation~\cite{QB-KS-GC:15} networks, to name only a few.  By
locally expressible, we mean that, for each
constraint, e.g., $g_j$, there exists an agent, which we
term corresponding agent, such that the function $g_j$
depends on the state of the corresponding agent and its
1-hop neighbors' state.  We assume that all the agents involved in a
constraint know the functional form of the constraint and its
derivatives.  According to this definition, different
  constraints might have different corresponding agents.  Under this
structure, agents require up to 2-hop communication to evaluate any
constraint in which they are involved (1-hop communication in the case
of the corresponding agent, 2-hop communication in the case of the
other agents involved in the constraint).

Our aim is to develop a smooth distributed algorithm to find an
optimizer of the constrained problem~\eqref{eq:nlp}.  Our solution
strategy employs a continuously differentiable exact penalty function,
cf. Section~\ref{sec:prelim}, to reformulate the problem as an
unconstrained optimization one. We then face the task of implementing
its gradient dynamics in a distributed way. To do so, we show that the
problem of distributed calculation of Lagrange multiplier functions
and other necessary terms in the gradient of the penalty function can
be formulated as a linear algebraic equation with separable data (cf.
Section~\ref{sec:dist_formulation}). In turn, we justify how this
algebraic equation can be solved in a distributed manner (cf.
Section~\ref{sec:dist_linear}). Finally, we combine both sets of
results to propose a distributed algorithmic solution based on smooth
gradient descent to solve~\eqref{eq:nlp}.

\begin{remark}\longthmtitle{Alternative approaches}
   To solve problem~\eqref{eq:nlp} in a distributed way, we can
    instead construct the Lagrangian and then use primal-dual (also
    known as saddle-point)
    dynamics~\cite{KA-LH-HU:58,DF-FP:10,AC-BG-JC:17-sicon}. This
    dynamics uses gradient descent in the primal variable and gradient
    ascent in the dual variable. For the problem structure described
    above, these dynamics is distributed (requiring up to 2-hop
    communication). However, the dynamics is in general slow, exhibits
    oscillations in the distance from the feasible set, and there is
    no guarantee of satisfying the constraints during the evolution,
    even if the initial state is feasible.  Also, it is not clear how
    to apply accelerated methods, cf.~\cite{YEN:83} to the primal-dual
    approach. Another approach to solve~\eqref{eq:nlp} in an (up to
    2-hop) distributed way consists of reformulating the problem as an
    unconstrained optimization~\cite{LX-SB:06,BJ-MJ:09,AC-JC:15-tcns}
    by adding to the original objective function non-differentiable
    penalty terms replacing the constraints~\cite{DPB:82} and
    employing subgradient-based methods. However, these methods are
    difficult to implement, often lead to chattering, and the study of
    their convergence properties requires tools from nonsmooth
    analysis.  Yet another approach is the alternating direction
    method of multipliers~\cite{SB-NP-EC-BP-JE:11}, which requires
    using some additional reformulation
    techniques~\cite{JFCM-JMFX-PMQA-MP:13} to make it distributed and
    convergence to an optimizer is only guaranteed when the
    optimization problem is convex.  Although it enjoys fast
    convergence, each agent needs to solve a local optimization
    problem at every iteration to update its state, which might be
    computationally inefficient depending on the form of the
    constraint and the objective functions.  \oprocend
\end{remark}

\section{Linear Algebraic Equations Defined by Separable Problem Data}\label{sec:dist_linear}
In this section, we propose a novel exponentially fast distributed
algorithm to solve linear algebraic equations whose problem data is
separable. As we argue later, such linear equations arise naturally
when considering the distributed solution of exact penalty
optimization problems, but the discussion here is of independent
interest.

Given a group of agents, consider a system of linear equations whose
coefficient matrix and constant vector are the aggregation of
individual coefficient matrices and constant vectors, one per
agent. Formally,
\begin{align}\label{eq:linear}
  \Big( \sum\limits_{i=1}^n N_i \Big) \, v &= \sum\limits_{i=1}^n b_i ,
\end{align}
where $n$ is the number of agents, $v \in \real^q$ is the unknown
solution vector, and $N_i \in \real^{q \times q}$ and
$b_i \in \real^q$ are the coefficient matrix and constant vector
corresponding to agent~$i$.
Our approach is based on first reformulating~\eqref{eq:linear} into a
1-hop distributed system of equations. By 1-hop (respectively, 2-hop)
distributed, we mean that each equation in the system only involves
some corresponding agent and its neighbors (respectively, 2-hop
neighbors).

We start by endowing each agent with its own candidate version
$v_i \in \real^q$ of $v$. Then, assuming a connected communication
graph among the agents, equation~\eqref{eq:linear} can be equivalently
rewritten as
\begin{subequations}\label{eq:realgebraic}
  \begin{align}
   \sum\limits_{i=1}^n N_i v_i
    &= \sum\limits_{i=1}^n
      b_i \label{eq:re1}
    \\
    (L \otimes I_{q}) \vv &=\mathbf{0}, \label{eq:re2}
  \end{align}
\end{subequations}
where $\vv=[v_1;\ldots;v_n] \in \real^{nq}$.  Note that, in order
for~\eqref{eq:re2} to be true, it must hold that
$\vv = \mathbf{1} \otimes v$, i.e., all $v_i$'s are the same.
Although~\eqref{eq:re2} is 1-hop distributed, \eqref{eq:re1} is not
distributed.  To address this, we introduce a new variable
$y_i \in \real^{q}$ per agent $i \in \until{n}$. Let
$\yy=[y_1 ;\ldots; y_n] \in \real^{nq}$ and consider the following set
of equations,
\begin{align}\label{eq:matrix}
  \underbrace{\begin{pmatrix}
      \N & -L \otimes I_{q}
      \\
      L \otimes I_{q} & \mathbf{0}
  \end{pmatrix}
  }_{\PP}
  \begin{pmatrix}
    \vv
    \\
    \yy
  \end{pmatrix}
  = 
  \underbrace{
  \begin{pmatrix}
    \bb
    \\
    \mathbf{0}
  \end{pmatrix}
  }_{\qq},
\end{align}
where $\N=\diag(N_i) \in \real^{nq \times nq}$ and $\bb=[ b_1 ; \ldots
; b_n] \in \real^{nq}$.  Note that the set of
equations~\eqref{eq:matrix} is 1-hop distributed. The
  following result characterizes the equivalence
  between~\eqref{eq:matrix} and~\eqref{eq:linear}.

\begin{proposition}\longthmtitle{Equivalence between~\eqref{eq:matrix}
    and~\eqref{eq:linear}}\label{prop:equi}
  The solutions of~\eqref{eq:matrix} are of the form
  $(\mathbf{1} \otimes v; \bar{\yy} + \mathbf{1} \otimes y)$, where
  $v \in \real^q$ solves~\eqref{eq:linear},
  $\N (\one \otimes v) - \bb = (L \otimes I_q) \bar{\yy}$, and
  $y \in \real^q$.
\end{proposition}
\begin{IEEEproof}
Note that~\eqref{eq:matrix} can be rewritten as
  \begin{subequations}
    \begin{align}
      \begin{bmatrix}
        N_1 & & \\
        & \ddots & \\
        & & N_n
      \end{bmatrix} \vv &=\begin{bmatrix} b_1 \\ \vdots \\ b_n
      \end{bmatrix}+  (L \otimes I_{q}) \yy   , \label{eq:mat1}
      \\
      (L \otimes I_{q}) \vv &=\zero. \label{eq:mat2}
    \end{align}
  \end{subequations}
  Equation~\eqref{eq:mat2} implies that $\vv=\one \otimes v$, with $v
  \in \real^m$.  Then, from~\eqref{eq:mat1}, we have for each $i \in
  \until{n}$,
  \begin{align*}
    N_i v = b_i + (L_i \otimes I_q) \yy ,
  \end{align*}
  where $L_i$ denotes the $i$th row of the Laplacian~$L$.  Summing
  over all agents, we obtain
  \begin{align*}
    \Big( \sum\limits_{i=1}^n N_i \Big) v = \sum\limits_{i=1}^n b_i+
    \sum\limits_{i=1}^n(L_i \otimes I_q) \yy.
  \end{align*}
  Since $\mathbf{1}'L=0$, the last summand vanishes, which
  yields~\eqref{eq:linear}. The expression for $\yy$ again follows
  directly from the fact that $\one'L=0$.
 \end{IEEEproof}

Our next goal is to synthesize a distributed algorithm to
solve~\eqref{eq:matrix}. Our algorithm design is based on formulating
this equation as an unconstrained optimization problem.  Let $\zz=(\vv;\yy)$
and consider the quadratic function $V_1:\real^{2qn}\rightarrow\real$
\begin{align}\label{eq:V}
  V_1(\zz)&=\dfrac{1}{2}(\PP \zz - \qq)' (\PP \zz - \qq).
\end{align}
Note that $V_1$ vanishes over the solution set of $\PP \zz = \qq$ and
takes positive values otherwise.  The problem of
solving~\eqref{eq:matrix} can be reformulated as
\begin{align*}
  \min_{\zz} V_1(\zz) .
\end{align*}
The gradient descent dynamics of $V_1$ is given by
\begin{align*}
  \dot{\zz}=-\PP'(\PP \zz - \qq) .
\end{align*}
When convenient, we refer to this dynamics as $\linear$.  In expanded
form, it takes the form
\begin{subequations}\label{eq:algo}
  \begin{align}
    \dot{\vv} &= -\N'[\N \vv \!-\!(L \otimes I_{q})\yy \!-\! \bb] \!-\!
              (L^2 \otimes I_{q}) \vv
    \\
    \dot{\yy} &= (L \otimes I_{q})[\N \vv -(L \otimes I_{q})\yy- \bb] .
  \end{align}
\end{subequations}
From~\eqref{eq:algo}, each agent $i \in \until{n}$ has the
  dynamics
  \begin{align*}
    \dot{v}_i &= -N_i'\Big(N_i v_i - b_i - \sum\limits_{j \in \NN_i}
    (y_i - y_j)\Big) - \sum\limits_{j \in \NN_i} (v_i^L-v_j^L)
    \\
    \dot{y}_i &= \sum\limits_{j \in \NN_i} \Big( N_i v_i - b_i -(N_j
    v_j - b_j)\Big) - \sum\limits_{j \in \NN_i} (y_i^L - y_j^L),
  \end{align*}
  where $v_k^L= \sum\limits_{j \in \NN_k} (v_k - v_j)$ and $y_k^L=
  \sum\limits_{j \in \NN_k} (y_k - y_j)$.  
    This algorithm is 2-hop distributed, meaning that to execute it,
  each agent~$i \in \until{n}$ needs to know its state $(v_i;y_i)$ and
  the state of its 2-hop neighbors. The next result characterizes its
convergence properties.

\begin{proposition}\longthmtitle{Exponential convergence
    of~\eqref{eq:algo} to solution of linear system}\label{prop:expo}
  The dynamics~\eqref{eq:algo} converges to a solution
  of~\eqref{eq:matrix} exponentially with a rate proportional
  to~$\lambda_2(\PP' \PP)$.
\end{proposition}
\begin{IEEEproof}
  Let $\ww \in \nulls(\PP)$ and note
  \begin{align*}
    \ww' \dot{\zz}=-\ww '\PP'(\PP \zz - \qq)=0 .
  \end{align*}
  This means that the dynamics of $\zz$ is orthogonal to
  $\nulls(\PP)$.  Let us decompose $\zz(t)$ as
  $\zz(t)=\zz_{\|}(t)+\zz_{\perp}(t)$.  Here, $\zz_{\|}(t)$ is the
  component of $\zz(t)$ in $\nulls(\PP)$ and $\zz_{\perp}(t)$ is the
  component orthogonal to it.  From the above discussion, we have that
  $\zz_{\|}(t)=\zz_{\|}(0)$ under the dynamics~\eqref{eq:algo}.  Since
  this component does not change, consider the particular solution
  $\zz^*$ of~\eqref{eq:matrix} that satisfies
  $\zz^*_{\|} = \zz_{\|}(0)$. Note that $\zz^*$ defined in this way is
  unique.  Now, consider the Lyapunov function
  $V_2: \real^{2qn} \rightarrow \real$
  \begin{align}\label{eq:V2}
    V_2(\zz)=\dfrac{1}{2}(\zz-\zz^*)'(\zz-\zz^*).
  \end{align}
  The derivative of $V_2$ along the dynamics~\eqref{eq:algo} is given by
  \begin{align*}
    L_{\linear}V_2(\zz)
    &=(\zz-\zz^*)'\dot{\zz}
    \\
    &=-(\zz-\zz^*)'\PP'(\PP \zz - \qq)
    \\
    &=-(\zz-\zz^*)'\PP'\PP(\zz-\zz^*) \leq -2\lambda_2 (\PP' \PP)V_2(\zz).
  \end{align*}
  The last inequality follows from the fact that the evolution of
  $\zz$ is orthogonal to the nullspace of $\PP$. This proves that,
  starting from $\zz(0)$, the dynamics converges to the solution
  $\zz^*$ of~\eqref{eq:matrix} exponentially fast with a rate
  determined by the minimum non-zero eigenvalue of~$\PP' \PP$.  
\end{IEEEproof}

Next, we examine the robustness to disturbances of the
dynamics~\eqref{eq:algo}. This is motivated by the observation that,
in practical scenarios, one may face errors in the execution due to
imperfect knowledge of the problem data, imperfect information about
the state of other agents, or other external disturbances.  Formally,
we consider
\begin{align}\label{eq:disturb}
  \dot{\zz}=\linear(\zz)+d(t)=-\PP'(\PP \zz - \qq)+d(t) ,
\end{align}
where $d(t)$ denotes the disturbance.

\begin{proposition}\longthmtitle{Robustness of~\eqref{eq:disturb} against
    disturbances}\label{prop:robust}
  The dynamics~\eqref{eq:disturb} is input-to-state stable (ISS) with
  respect to the set of equilibria of~\eqref{eq:algo}.
\end{proposition}
\begin{IEEEproof}
  The disturbance $d(t)$ in~\eqref{eq:disturb} can be decomposed as
  $d(t) = d_{\|}(t)+d_{\perp}(t)$.  Due to the presence of
  $d_{\|}(t) \in \nulls(\PP)$, the component of $\zz(t)$ in
  $\nulls(\PP)$ does not remain constant any more. In fact,
  along~\eqref{eq:disturb}, we have $\ww' \dot \zz = \ww' d_{\|}$ for
  all $\ww \in \nulls(\PP)$, and therefore we deduce that
  $\dot \zz_{\|}(t) = d_{\|}(t)$.
   Consider then the equilibrium trajectory $t\mapsto z^*(t)$, where
  $\zz^*(t)$ is uniquely determined by the equations
  $\PP \zz^*(t) = \qq$ and $\zz^*_{\|}(t) = \zz_{\|}(t)$.
  Let $V_2$ be the same function as in~\eqref{eq:V2}, but now with the
  time-varying $\zz^*(t)$. The derivative of $V_2$ is given by
  \begin{align*}
    L_{\linear +d}V_2
    &=
      (\zz-\zz^*)'(\dot{\zz}-\dot{\zz}^*)
    \\ 
    & =  (\zz-\zz^*)'(-\PP'(\PP \zz - \qq)+d-d_{\|})
    \\ 
    & =  (\zz-\zz^*)'(-\PP'(\PP \zz - \qq)+d_{\perp})
       \\
    & \leq  -\lambda_2 (\PP' \PP)\|\zz-\zz^*\|^2+\|\zz-\zz^*\|\|d_{\perp}\|
    \\
    & \leq  -\lambda_2 (\PP' \PP)\|\zz-\zz^*\|^2+\|\zz-\zz^*\|\|d\|.
  \end{align*}    
  Choose $\theta \in (0,1)$. Then the above inequality can be
  decomposed as
  \begin{align*}
    L_{\linear+d}V_2
    & \leq
      -\lambda_2 (\PP' \PP)(1-\theta)\|\zz-\zz^*\|^2 
    \\
    & \quad 
      - \lambda_2 (\PP' \PP)\theta\|\zz-\zz^*\|^2 +\|\zz-\zz^*\|\|d\|.
  \end{align*}
  Hence, $L_{\linear+d}V_2\leq -\lambda_2 (\PP' \PP)(1-\theta)\|\zz-
  \zz^*\|^2$ if $\|\zz- \zz^*\| \geq \dfrac{\|d\|}{\lambda_2 (\PP'
    \PP)\theta}$.  From~\cite[Theorem 4.19]{HKK:02}, this means that
  the system is input-to-state stable with respect to the set of
  equilibria with gain $\gamma(r)=\dfrac{r}{\lambda_2 (\PP'
    \PP)\theta}$.  
\end{IEEEproof}

Proposition~\ref{prop:robust} implies that the trajectories
of~\eqref{eq:disturb} asymptotically converge to a neighborhood of the
set of equilibria of~\eqref{eq:algo} (with the size of the
neighborhood scaling up with the size of the disturbance).  All
equilibria correspond to solutions of~\eqref{eq:linear}.  The results
of this section show that the system~\eqref{eq:linear} can be solved
in a distributed and robust way.

  \begin{remark}\longthmtitle{Distributed algorithms for linear
      algebraic equations}  Although we consider the linear
      algebraic equations~\eqref{eq:linear} here to perform the
      distributed computation of the gradient of the penalty function,
      solving linear algebraic equations in a distributed fashion is
      an interesting problem on its own,
      cf.~\cite{DPB-JNT:97,SM-JL-ASM:15,BDOA-SM-ASM-UH:16}. Different
      algorithmic solutions exist depending on the assumptions about
      the information available to the individual agents.
      Specifically, equations with the same structure
      as~\eqref{eq:linear} appear frequently~\cite{JL-CYT:16} with
      applications to distributed sensor
      fusion~\cite{DPS-ROS-RMM:05-ifac} and maximum-likelihood
      estimation~\cite{LX-SB-SL:05}. \cite{DPS-ROS-RMM:05-ifac}
      exploits the positive definiteness of the matrices and
      \cite{LX-SB-SL:05} uses element-wise average consensus to find
      the solutions of~\eqref{eq:linear}. \cite{JL-CYT:16} also
      exploits the positive definite property of the individual
      matrices and requires the agents to know the state as well as
      the matrices of the neighbors.  The algorithmic design procedure
      we employ here is similar to the one used in~\cite{XW-SM:18},
      which leads to an algorithm that also does not require the
      positive definiteness of the individual matrices.
      Interestingly, the convergence analysis in~\cite{XW-SM:18} uses
      the linearity of the dynamics and La Salle's invariance
      principle to conclude exponential stability, although it does
      not guarantee that the agents converge to the same solution.  By
      contrast, the Lyapunov-based technical analysis presented here,
      based on exploiting the orthogonality of the dynamics to the
      nullspace of the reformulated system matrix, allows us to lower
      bound the exponential convergence rate and formally characterize
      the robustness properties of the algorithm against
      disturbances. Both properties are key for the application later
      in Section~\ref{sec:dist_optimization} to distributed gradient
      computation via characterizing the stability of the
      interconnected system.   \oprocend
\end{remark}

\section{Distributed Computation of the Gradient of Penalty
  Function}\label{sec:dist_formulation}

We pursue next our strategy to solve the constrained optimization
problem~\eqref{eq:nlp} in a distributed fashion by using the gradient
dynamics of the continuously differentiable exact penalty
function~\eqref{eq:penalty}.  In this section, we first identify the
challenges associated with the distributed computation of
$\nabla f^\epsilon$ and then employ the algorithmic tools and results
of Section~\ref{sec:dist_linear} to address them.

The gradient of $f^{\epsilon}(x)$ with respect to $x_i$ is given by
\begin{align}\label{eq:grad_single}
  \nabla_{x_i} f^\epsilon(x)
  &=\nabla_{x_i}f_i(x_i)+\sum\limits_{j=1}^m
    \lambda_j(x) \nabla_{x_i} g_j(x)
  \\
  & \quad + \sum\limits_{k=1}^p \mu_k(x)
    \nabla_{x_i} h_k(x) + \sum\limits_{k=1}^p h_k(x) \nabla_{x_i}
    \mu_k(x) \notag
  \\
  & \quad + \sum\limits_{j=1}^m
    \big(g_j(x)+y_j^{\epsilon2}(x)\big) \nabla_{x_i} \lambda_j(x) 
    \notag
  \\
  & \quad + \frac{2}{\epsilon} \sum\limits_{j=1}^m
    \big(g_j(x)+y_j^{\epsilon2}(x)\big) \nabla_{x_i} g_j(x) \notag
  \\
  & \quad
    + \frac{2}{\epsilon} \sum\limits_{k=1}^p h_k(x) \nabla_{x_i}
    h_k(x) .
    \notag
\end{align}
In this expression, and with the assumptions made in
Section~\ref{sec:problem}, if agent $i$ knew $(\lambda(x),\mu(x))$,
then it could compute all the terms locally except for
\begin{align*}
  \rho_i (x) \equiv \sum\limits_{j=1}^m
  \big(g_j(x)+y_j^{\epsilon2}(x)\big) \nabla_{x_i}
  \lambda_j(x)+\sum\limits_{k=1}^p h_k(x) \nabla_{x_i} \mu_k(x) .  
\end{align*}
The rest of this section is devoted to show how to deal with these two
issues. First, we show how we can formulate and solve the problem of
calculating $(\lambda(x),\mu(x))$ in a distributed way. After that, we
show how agent $i$ can calculate $\rho_i(x)$ with only local
information and communication.

\subsection{Distributed computation of multiplier functions}

Given $x \in \real^n$, $(\lambda(x),\mu(x))$ are defined by the linear
algebraic equation~\eqref{eq:lambda}. Note that this equation can also be
written as
\begin{align}
  N(x)
  \begin{bmatrix}
    \lambda(x)
    \\
    \mu(x)
  \end{bmatrix}
  &= - 
  \begin{bmatrix}
    \nabla g(x)'
    \\
    \nabla h(x)'
  \end{bmatrix}
  \nabla f(x) . \label{eq:relambda}
\end{align}
The next result proves that we can actually decompose the matrix
$N(x)$ and the righthand side of~\eqref{eq:relambda} as the summation
of locally computable matrices. This makes the equation have the same
structure as equation~\eqref{eq:linear}, and hence we can use the distributed
algorithm in Section~\ref{sec:dist_linear} to solve it.

\begin{proposition}\longthmtitle{Equivalence
    between~\eqref{eq:relambda}
    and~\eqref{eq:linear}}\label{prop:lambdamu}
  For each $x \in \real^n$, calculating $(\lambda(x),\mu(x))$ can be
  cast as solving a linear algebraic equation of the
  form~\eqref{eq:linear}.
\end{proposition}
\begin{IEEEproof}
  For convenience, for each $i \in \until{n}$, we define
  $\mathrm{g}_i(x)=(\mathrm{g}_{i1}(x),\dots,\mathrm{g}_{im}(x)) \in
  \real^m$, where
  \begin{align}\label{eq:g}
    \mathrm{g}_{ij}(x)=
    \begin{cases}
      \frac{g_j(x)}{n_j} & \text{if $i$ is involved in constraint $g_j$}
      \\
      0 & \text{otherwise}
    \end{cases}
  \end{align}
  and $n_j$ is the total number of agents involved in constraint $j
  \in \until{m}$. From this definition, we have that $\sum_{i=1}^n
  (\diag( \gamma \mathrm{g}_i(x)))^2 = \gamma^2 G^2(x)$. Using this
  fact, we define $N_i(x)$, for each $i \in \until{n}$ as
  \begin{align*}
    N_i(x) \!\!=\! \!
    \begin{bmatrix}
      \nabla_{x_i}g(x)'
      \\
      \nabla_{x_i}h(x)'
    \end{bmatrix}
    \!\!
    \begin{bmatrix}
      \nabla_{x_i}g(x) \; \nabla_{x_i}h(x)
    \end{bmatrix} \!+\! 
    \begin{bmatrix}
      (\diag( \gamma \mathrm{g}_i(x)))^2 &\!\!\! \mathbf{0}
      \\
      \mathbf{0} & \!\!\!\mathbf{0}
    \end{bmatrix}
  \end{align*}
  From the definition~\eqref{eq:n} of $N(x)$, note that
  \begin{align}\label{eq:distN}
    N(x)=\sum\limits_{i=1}^n N_i(x) .
  \end{align}

  The righthand side of~\eqref{eq:relambda} could be decomposed as
  \begin{align}\label{eq:distB}
    \begin{bmatrix}
      \nabla g(x)' \\ \nabla h(x)'
    \end{bmatrix} \nabla f(x)=\sum\limits_{i=1}^n \begin{bmatrix}
    \nabla_{x_i}g(x)' \nabla_{x_i} f_i(x_i) \\
    \nabla_{x_i}h(x)' \nabla_{x_i} f_i(x_i)
  \end{bmatrix}.
  \end{align}
  Hence,~\eqref{eq:relambda} is equivalent to \eqref{eq:linear} with
  $q=m+p$, completing the proof. 
\end{IEEEproof}

From the definition of the matrices $\{N_i\}_{i=1}^n$ in the proof of
Proposition~\ref{prop:lambdamu}, one can deduce that, individually,
these matrices might not be positive definite in general. This
highlights the importance of the algorithm~\eqref{eq:algo} to solve
equations with separable problem data when the coefficient matrices
are not necessarily positive definite.  Combining
Proposition~\ref{prop:lambdamu} with the discussion of
Section~\ref{sec:dist_linear}, we deduce that each agent can compute
$(\lambda(x), \mu(x))$ in a distributed~way.

\subsection{Distributed computation of the gradient}

Here, we describe how agent $i \in \until{n}$ can calculate
$\rho_i(x)$ locally, completing the distributed computation of
$\nabla_{x_i}f^\epsilon(x)$.

\begin{proposition}\longthmtitle{Local computation of
    $\rho_i(x)$}\label{prop:rho}
  For each $x\in \real^n$, agent $i \in \until{n}$ can calculate
  $\rho_i(x)$ locally via communication with its 2-hop neighbors.
\end{proposition}
\begin{IEEEproof}
  In compact form, $\rho_i(x)$ can be written as
  \begin{align*}
    \rho_i(x)=
    \begin{bmatrix}
      g(x) + Y^\epsilon(x) y^\epsilon(x)
      \\
      h(x)
    \end{bmatrix}'
    \begin{bmatrix}
      \nabla_{x_i} \lambda(x)
      \\
      \nabla_{x_i} \mu(x)
    \end{bmatrix} .
  \end{align*} 
  This means that $\rho_i(x)$ is given by the $i$th column of $[ g(x)
  + y^\epsilon(x); h(x)]' [ \nabla \lambda(x)' ; \nabla
  \mu(x)']$. From~\eqref{eq:grad_lambda}, this is equivalent to saying
  that $\rho_i(x)$ is given by $- \varrho(x)' ( r_i(x); s_i(x))$,
  where $\varrho(x)=N^{-1}(x) [ g(x) +Y^\epsilon(x) y^\epsilon(x);
  h(x)]$ whose transpose is $[ g(x) +Y^\epsilon(x) y^\epsilon(x);
  h(x)]' N^{-1}(x)$ (since $N(x)$ is symmetric) and $( r_i(x);
  s_i(x))$ denotes the $i$th column of $[R(x); S(x)]$.
  Based on this, we divide the distributed computation of $\rho_i(x)$
  in two parts:
  \begin{enumerate}[(a)]
  \item First we show how all agents can compute $\varrho(x)$ using a 2-hop distributed
    algorithm;
  \item Next we show that each agent $i \in \until{n}$ can calculate
    $r_i(x)$ and $s_i(x)$ locally via communication with its 2-hop
    neighbors.
  \end{enumerate}
  
  For (a), consider the following equation in $\varrho$
  \begin{align}\label{eq:rho}
    N(x) \varrho(x) = \begin{bmatrix}
      g(x) +Y^\epsilon(x) y^\epsilon(x)\\
      h(x)
    \end{bmatrix}.
  \end{align}
  We can decompose the righthand side
  of~\eqref{eq:rho}~as
  \begin{align}\label{eq:distgh}
    \begin{bmatrix}
      g(x)+Y^\epsilon(x)y^\epsilon(x)
      \\
      h(x)
    \end{bmatrix}
    = \sum\limits_{i=1}^n
    \begin{bmatrix}
      \mathrm{g}_i(x) + \mathrm{y}^2_i(x)
      \\
      \mathrm{h}_i(x)
    \end{bmatrix} ,
  \end{align}
  where $\mathrm{g}_i(x)$ is defined in~\eqref{eq:g}, and
  $\mathrm{y}^2_i(x)$ and $\mathrm{h}_i(x)$ are defined
  similarly. From~\eqref{eq:distN} and~\eqref{eq:distgh},
  equation~\eqref{eq:rho} has the structure described
  in~\eqref{eq:linear} and hence can be solved in a distributed manner
  by the algorithm of Section~\ref{sec:dist_linear}.

  Next we look at the decomposition of $[R(x) ; S(x)]$ for (b). We
  describe here only the decomposition for $R(x)$ (the decomposition
  for $S(x)$ is similar). From~\eqref{eq:r}, $R(x)$ in
  expanded form is
  \begin{align*}
    &\nabla g(x)' \Big(\nabla^2 f(x) + \sum\limits_{j=1}^m
      \lambda_j(x) \nabla^2 g_j(x) + \sum\limits_{k=1}^p \mu_k
      \nabla^2 h_k(x) \Big)
    \\
    &+\sum\limits_{j=1}^m e_j^m (\nabla f(x)' + \lambda' \nabla g(x)'
    + \mu' \nabla h(x)') \nabla^2 g_j(x)
    \\
    & + 2\gamma^2 \Lambda(x) G(x) \nabla g(x)',
  \end{align*}
  which clearly corresponds to a sum of matrices. Here, we look at the
  first column of these matrices one by one and show that $r_1(x)$ can
  be calculated by agent 1 with information from its 2-hop neighbors
  (following the same reasoning justifies that each $r_i(x)$ can be
  calculated by agent $i \in \until{n}$).  The first column of the
  first matrix is given by $\nabla g(x)' \nabla^2_{x_1 x} f(x)$. To
  calculate it, in addition to $\nabla^2_{x_1} f_1(x)$, agent 1 only
  needs to know the partial derivative of the constraints in which it
  is involved (which are available to it by assumption,
  cf. Section~\ref{sec:problem}). The first column corresponding to
  the next two matrices is given by $\nabla g(x)' \Big(
  \sum\limits_{j=1}^m (\nabla^2_{x_1x}g_j(x))\lambda_j(x) +
  \sum\limits_{k=1}^p (\nabla^2_{x_1x}h_p(x))\mu_k(x) \Big)$.  For
  these, agent 1 only needs information about the partial first and
  second derivatives of the constraints in which it is involved, in
  addition to the values of the multiplier functions.
  The first column corresponding to the next three matrices is
  $\sum\limits_{j=1}^m e_j^m (\nabla f(x)' + \lambda' \nabla g(x)' +
  \mu' \nabla h(x)') \nabla^2_{x x_1} g_j(x). $ The calculation of the
  first term is straightforward. Rewriting the second term as
  $\sum\limits_{j=1}^m e_j^m \lambda' \nabla g(x)' \nabla^2_{x x_1}
  g_j(x)$ and knowing the structure of $\nabla g(x)' \nabla^2_{x x_1}
  g_j(x)$ from the discussion above, we can say that can be calculated
  by agent 1 (a similar observation applies to the third
  term). Regarding the last matrix, the first column is $2\gamma^2 [
  \lambda_1 g_1 \nabla_{x_1}g_1 ; \ldots ; \lambda_m g_m
  \nabla_{x_1}g_m]$. Clearly, agent 1 only needs to know the values
  and partial derivatives of the constraints in which it is involved
  for calculating this, concluding the proof.  
\end{IEEEproof}

\begin{remark}\longthmtitle{Scalability with the number of agents}
  In the preliminary conference version~\cite{PS-JC:18-cdc} of
    this work, we had all agents compute the Jacobian matrix of the
    multiplier functions to calculate $\rho_i(x)$. Since the dimension
    of the Jacobian matrix is $(m+p) \times n$, this approach was not
    scalable with the number of agents. With the approach described
    here, instead, each agent only needs to compute a vector of size
    $(m+p) \times 1$, which scales independently with the number of
    agents~$n$.  \oprocend
\end{remark}

Based on Propositions~\ref{prop:lambdamu} and~\ref{prop:rho}, for a
given $x \in \real^n$, we can compute asymptotically the values of
$\lambda(x)$, $\mu(x)$ and $\varrho(x)$, and in turn, the gradient of
the penalty function in a distributed way. For its use later, we
denote by $\PPest(x)$ the corresponding matrix defined as
in~\eqref{eq:matrix}, which now depends on $x$ due to the
$x$-dependence of $N_i$ and $b_i$ (and hence $\N$ and $\bb$) in
equations~\eqref{eq:relambda} and~\eqref{eq:rho}.

\begin{remark}\longthmtitle{Robustness in the calculation of gradient}\label{remark:robustness}
   From Proposition~\ref{prop:robust}, the distributed calculation
    of the gradient of the exact penalty function is robust to bounded
    disturbances due to errors in the problem data (e.g., errors in
    the value of the constraint functions or the gradients of the
    objective and constraint functions), packet drops, or
    communication noise. Furthermore, since the matrix
    $N(x)=\sum_{i}^n N_i(x)$ is positive definite (and hence
    invertible) from Proposition~\ref{prop:lambda}, it follows that
    all equilibria have the same unique variable $\vv$, whereas the
    auxiliary ones $\yy$ may take multiple values according
      to Proposition~\ref{prop:equi}.  This means that, for a given
    $x \in \real^n$, the primary variables $\{v_i\}_{i=1}^n$ converge
    uniquely to $\lambda(x)$, $\mu(x)$ and $\varrho(x)$ under each of
    the algorithms described above.   \oprocend
\end{remark}

\section{Distributed Optimization via Interconnected Dynamics}\label{sec:dist_optimization}

In this section, we finally put all the elements developed so far
together to propose a distributed algorithm to solve~\eqref{eq:nlp}.
The basic idea is to implement the gradient dynamics of the exact
penalty function. However, the algorithmic solutions resulting from
Section~\ref{sec:dist_formulation} only asymptotically compute the
gradient of the exact penalty function at a given state. This state,
in turn, changes by the action of the gradient descent dynamics. The
proposed distributed algorithm is then the result of the
interconnection of these two complementary dynamics.

Formally, the gradient descent dynamics of $f^\epsilon$ which serves
as reference for our algorithm design takes the form
\begin{align}\label{eq:idealgd}
  \dot x = -\nabla f^\epsilon(x) .
\end{align}
For convenience, define $\chi: \real^n \to \real^{2(m+p)}$ by $\chi(x) = (\lambda(x), \mu(x),\varrho(x))$ and rewrite~\eqref{eq:idealgd} as $
\dot{x}=\grad(x,\chi(x))$ for an appropriate function $\grad$ defined
by examining the expression in~\eqref{eq:grad_single} for $i \in
\until{n}$ (note that, given the assumptions on the problem functions,
for each $x \in \real^n$, the function $\grad$ is locally Lipschitz in
its argument~$\chi$, and from Proposition~\ref{prop:lambda} and
equation~\eqref{eq:rho}, $x \mapsto \chi(x)$ is continuously
differentiable).  The variable $\chi$ corresponds to those terms
appearing in the gradient that are not immediately computable with
local information.  However, with the distributed algorithms described
in Section~\ref{sec:dist_formulation}, the network agents can
asymptotically compute $\chi(x)$ in a distributed fashion.  Let
$\Upsilon \in \real^{4(m+p)}$ denote the augmented variable containing the estimates of
$\chi(x)$ and the associated auxiliary variables, available to the
network agents via
\begin{subequations}\label{eq:algorithm}
  \begin{align}\label{eq:estimation}
     \dot{\Upsilon}&=\est (x, \Upsilon) ,
  \end{align}
  where $\est(x,\Upsilon)$ denotes the algorithms of the
  form~\eqref{eq:algo} described in Section~\ref{sec:dist_formulation}. 
  Let $ \hat{\chi} = \mathcal{P}_{\chi} \Upsilon$ denote the
  projection of $\Upsilon$ onto the $\chi$ space, i.e., corresponding
  to the set of primary variables. From Proposition~\ref{prop:expo},
  we note that, for fixed $x\in \real^n$, $\hat{\chi} \to \chi(x)$
  exponentially fast.  Hence, with the information available to the
  agents, instead of~\eqref{eq:idealgd}, the network implements
  \begin{align}\label{eq:gdynamics}
    \dot{x}= \grad(x,\mathcal{P}_{\chi} \Upsilon) .
  \end{align}
\end{subequations}
Our proposed algorithm is the interconnected dynamical
system~\eqref{eq:algorithm}. When convenient, we refer to it as~$\ic$.
Note that this algorithm is 2-hop distributed. Moreover, for each
equilibrium $(x_{\text{eq}},\Upsilon_{\text{eq}})$,
of~\eqref{eq:algorithm}, its $x$-component $x_{\text{eq}}$ is an
equilibrium of~\eqref{eq:idealgd} (which is also a KKT point of
problem~\eqref{eq:nlp} if EMFCQ is satisfied, cf.
Proposition~\ref{prop:exactness}). We characterize the convergence
properties of the algorithm~\eqref{eq:algorithm} next.

\begin{theorem}\longthmtitle{Asymptotic convergence of distributed algorithm
    to solution of optimization problem}\label{thm:convergence-algo}
  Assume LICQ is satisfied at each $x \in \mathcal{D}$.  For each $x$,
  let $L_{\chi}(x)$ be the Lipschitz constant of $\chi \mapsto
  \grad(x,\chi)$.  Then the equilibria of the interconnected
  dynamics~\eqref{eq:algorithm} are asymptotically stable
  if there exists $\alpha > 0$ such that
  \begin{align}\label{eq:condition}
    \max_{x \in \D} \dfrac{\eta_{\alpha}(x)}{\lambda_2(\PPest(x)'
      \PPest(x))} < 1,
  \end{align}
  where $\eta_{\alpha}(x)=\frac{1}{4 \alpha} (\alpha L_{\chi}(x) + \|
  \nabla_x \chi (x) \| )^2 + L_{\chi}(x) \| \nabla_x \chi (x) \|$.
\end{theorem}
\begin{IEEEproof}
  We start by noting that~\eqref{eq:condition} is well defined since,
  from the definitions of $R(x)$ and $S(x)$ in~\eqref{eq:rs} and the
  expression of the gradient in~\eqref{eq:grad_single}, we deduce that
  $L_\chi(x)$ is continuous in $x$, and, moreover, since $\D$ is
  compact, $L_{\chi}$ and $\| \nabla_x \chi \|$ are bounded over $\D$.
  Consider now the Lyapunov function candidate for the interconnected
  system as
  \begin{align}\label{eq:Vc}
    V_{\text{c}} (x,\Upsilon) = \alpha f^\epsilon(x) + V_2(x,
    \Upsilon),
  \end{align}
  where $V_2$ is defined as in~\eqref{eq:V2}, but due to the
  dependence of $\zz^*$ on $x$ from equations~\eqref{eq:relambda}
  and~\eqref{eq:rho}, is now a function of $x$ too. The derivative
  of~$V$ along the dynamics~\eqref{eq:algorithm} is
  \begin{align*}
    & L_{\ic} V_{\text{c}} (x,\Upsilon) 
      \\
    & \quad = (\alpha \nabla f^\epsilon(x)+ \nabla_x V_2)'
    \grad(x,\mathcal{P}_{\chi} \Upsilon) + \nabla_{\Upsilon} V_2'
    \est(x, \Upsilon)
    \\
    & \quad \leq - (\alpha \nabla f^\epsilon(x)+ \nabla_x \chi
    (\hat{\chi} - \chi(x)) )' (\nabla f^\epsilon(x)
    \\
    & \qquad - \grad(x,\hat{\chi})+ \grad(x, \chi(x))) - \lambda_2 (x)
    \|\hat{\chi} - \chi(x) \|^2 ,
  \end{align*}
  where we have added and subtracted $\nabla f^\epsilon(x) = \grad(x,
  \chi(x))$ to $\grad(x,\mathcal{P}_{\chi} \Upsilon)$ and used the
  shorthand notation $\lambda_2(x) \equiv \lambda_2(\PPest'(x)
  \PPest(x))$.
  Hence, we have
  \begin{align*}
    L_{\ic} V_{\text{c}} (x,\Upsilon) \leq -\begin{bmatrix} \| \nabla
      f^\epsilon(x) \|
      \\
      \|\hat{\chi} - \chi(x) \|
    \end{bmatrix}'
    A(x)
    \begin{bmatrix} 
      \| \nabla f^\epsilon(x) \|
      \\
      \|\hat{\chi} - \chi(x) \|
    \end{bmatrix},
  \end{align*}
  with
  \begin{align*}
    A(x)=\begin{bmatrix}
      \alpha  & -\frac{1}{2} (\alpha L_{\chi}(x)+ \| \nabla_x \chi \|) \\
      -\frac{1}{2} (\alpha L_{\chi}(x)+ \| \nabla_x \chi \|) &
      \lambda_2 (x) - L_{\chi}(x) \| \nabla_x \chi \|)
    \end{bmatrix}.
  \end{align*}
  Next, we examine the positive-definiteness nature of the $2 \times
  2$-matrix $A(x)$. Since $\alpha>0$, note that $ A(x) \succ 0$ if the
  determinant is positive.  For $x \in \D$, the latter holds if and
  only if  $\alpha$ is such that
  \begin{align*}
    \eta_{\alpha}(x) / \lambda_2(x) < 1.
  \end{align*}
  Hence, under~\eqref{eq:condition}, this inequality holds over $\D$,
  and consequently $L_{\ic} V_{\text{c}} (x,\Upsilon) < 0$ over $\D
  \times \real^{4(m+p)}$.
 \end{IEEEproof}

The condition~\eqref{eq:condition} in
Theorem~\ref{thm:convergence-algo} can be interpreted as requiring the
estimation dynamics~\eqref{eq:estimation} to be fast enough to ensure
the error in the gradient computation remains manageable, resulting in
the convergence of the interconnected system.  In general,
however,~\eqref{eq:condition} might not be satisfied. To address this,
and inspired by this interpretation, we propose to execute the
estimation dynamics on a tunable timescale,
substituting~\eqref{eq:estimation} by
\begin{align}\label{eq:estimationtau}
  \tau \dot{\Upsilon}&=\est(x, \Upsilon) .
\end{align}
Here, $\tau>0$ is a design parameter capturing the timescale at which
the estimation dynamics is now executed.  Resorting to singular
perturbation theory, cf.~\cite{HKK:02,VV:97}, one could show that
$x(t) \to x_{\text{grad}}(t)$ as $\tau \to 0$, where $x_{\text{grad}}$
denotes the trajectory of the gradient descent
dynamics~\eqref{eq:idealgd}.  However, for the proposed approach to be
practical, it is desirable to have a strictly positive value of the
timescale below which convergence is guaranteed.  The following result
shows that such critical value exists.

\begin{proposition}\longthmtitle{Asymptotic convergence of distributed
    algorithm via accelerated estimation
    dynamics}\label{prop:convergence-algo}
  Assume LICQ is satisfied at each $x \in \mathcal{D}$ and let
  \begin{align*}
    \tau_* = \dfrac{ \lambda_{\min} (\PPest' \PPest)}{ 2
      \bar{L}_{\chi} \| \nabla_x \bar{\chi} \|} > 0,
  \end{align*}
  where $\lambda_{\min} (\PPest' \PPest)$ denotes the minimum of
  $\lambda_2 (\PPest'(x) \PPest(x))$, and $\bar{L}_{\chi}$ and $\|
  \nabla_x \bar{\chi} \|$ denote the maximum of $L_{\chi}$ and $\|
  \nabla_x \chi \|$ resp., over $\D$. Then, for any $\tau \in [0,
  \tau_*)$, the equilibria of the interconnected
  dynamics~\eqref{eq:gdynamics} and~\eqref{eq:estimationtau} are
  asymptotically stable.
\end{proposition}
\begin{IEEEproof}
  Let $\alpha>0$ and consider the Lyapunov function
  candidate~\eqref{eq:Vc}. Define
  \begin{align*}
    A_{\tau}(x)\!\! =\!\!
    \begin{bmatrix}
      \alpha & -\frac{1}{2} (\alpha L_{\chi}(x)+ \| \nabla_x \chi \|)
      \\
      -\frac{1}{2} (\alpha L_{\chi}(x)+ \| \nabla_x \chi \|) &
      \tau^{-1} \lambda_2 (x) - L_{\chi}(x) \| \nabla_x \chi \|)
    \end{bmatrix}\!\!.
  \end{align*}
  Following the same line of argument as in the proof of
  Theorem~\ref{thm:convergence-algo}, we arrive at
  \begin{align*}
    L_{\ic} V_{\text{c}} (x,\Upsilon) \leq -
    \begin{bmatrix}
      \| \nabla f^\epsilon(x) \|
      \\
      \|\hat{\chi} - \chi(x) \|
    \end{bmatrix}' A_{\tau}(x)
    \begin{bmatrix}
      \| \nabla f^\epsilon(x) \|
      \\
      \|\hat{\chi} - \chi(x) \|
    \end{bmatrix} 
  \end{align*}
  and the condition $\tau < \lambda_2(x) / \eta_{\alpha}(x) $ to
  ensure $L_{\ic}V_{\text{c}} (x,\Upsilon) < 0$.  Using the bounds for
  $L_{\chi}$ and $\| \nabla_x \chi \|$, we  upper bound
  $\eta_{\alpha} $ over $\D$ as 
  \begin{align*}
    \eta_{\alpha} (x) \le \bar{\eta}_{\alpha} = \frac{1}{4 \alpha}
    (\alpha \bar{L}_{\chi} + \| \nabla_x \bar{\chi} \| )^2 +
    \bar{L}_{\chi} \| \nabla_x \bar{\chi} \| .
  \end{align*}
  Consequently, it is enough to have $\tau < \lambda_2(x) /
  \bar{\eta}_{\alpha} $ for all $x \in \D$.  To establish the maximum
  admissible value of $\tau$, we can select the value of $\alpha$
  minimizing $\bar{\eta}_{\alpha}$. Since $\bar{\eta}_{\alpha}$ is
  strictly convex in $\alpha \in [0,\infty)$,
  this is given by the solution of
  \begin{align*}
    \dfrac{d}{d \alpha} \left( \dfrac{1}{ \alpha} (\alpha
      \bar{L}_{\chi} + \| \nabla_x \bar{\chi} \| )^2 \right) = 0.
  \end{align*}
  After some algebraic manipulations, one can verify that $\alpha^* =
  \| \nabla_x \bar{\chi} \| / \bar{L}_{\chi}$. Substituting this value
  in the expression of $\bar{\eta}_{\alpha}$ and taking the minimum
  over all $x \in \D$ yields  the definition of $\tau_*$.
\end{IEEEproof}

Note that the conditions identified in
Theorem~\ref{thm:convergence-algo} and
Proposition~\ref{prop:convergence-algo} to ensure convergence are
based on upper bounding the terms appearing in the Lie derivative of
the Lyapunov function candidate using 2-norms and, as such, are
conservative in general. In fact, the algorithm may converge even if
these conditions are not satisfied, something that we have observed in
simulation.

\begin{remark}\longthmtitle{Constraint satisfaction with the
    distributed dynamics}
  The centralized gradient descent on which we build our approach
  enjoys the constraint satisfaction properties stated in
  Proposition~\ref{prop:feasible}.  This means, using a singular
  perturbation argument~\cite{HKK:02,VV:97}, that the distributed
  gradient descent approach proposed here has the same guarantees as
  $\tau \to 0$.  Although for a fixed $\tau > 0$, we do not have a
  formal guarantee that the state remains feasible, we have observed
  this to be the case in simulations, even under general constraints.
  We believe this is due to the error-correcting terms in
    the original penalty function, which penalize deviations from the
    feasible set.  This anytime nature is especially important in
  applications where the optimization problem is not stand alone and
  its solution serves as an input to other layer in the control design
  (for example as a power/thermal set point, cf.~\cite{JR-HJ:14,AUR-SK:12}), where
  the algorithm should yield a feasible solution if terminated in
  finite time.  \oprocend
\end{remark}

\section{Simulations}\label{sec:sims}
\begin{figure}[htb]
  \centering
  \includegraphics[scale=0.45]{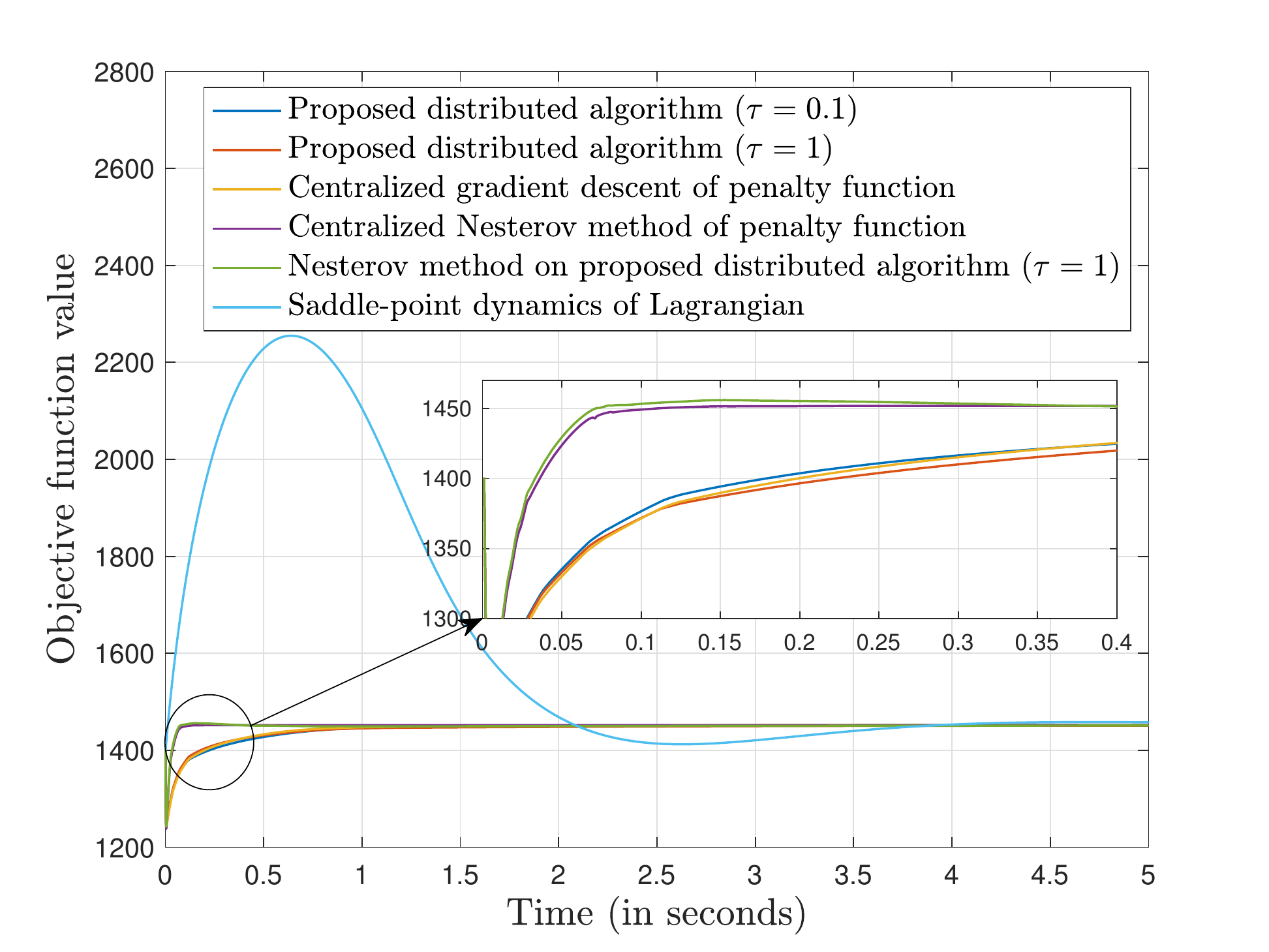}
  \caption{Evolution of the objective function value under the
    proposed distributed dynamics with $\tau=10^{-1}$ and $1$, resp.,
    the centralized gradient descent, the centralized and the
    distributed (using the proposed approach with $\tau=1$) Nesterov's
    accelerated gradient method of the penalty function, and the
    saddle-point dynamics of the Lagrangian.}\label{fig:comparison}
\end{figure}

\begin{figure}[htb]
  \centering
  \subfigure[]{\includegraphics[scale=0.45]{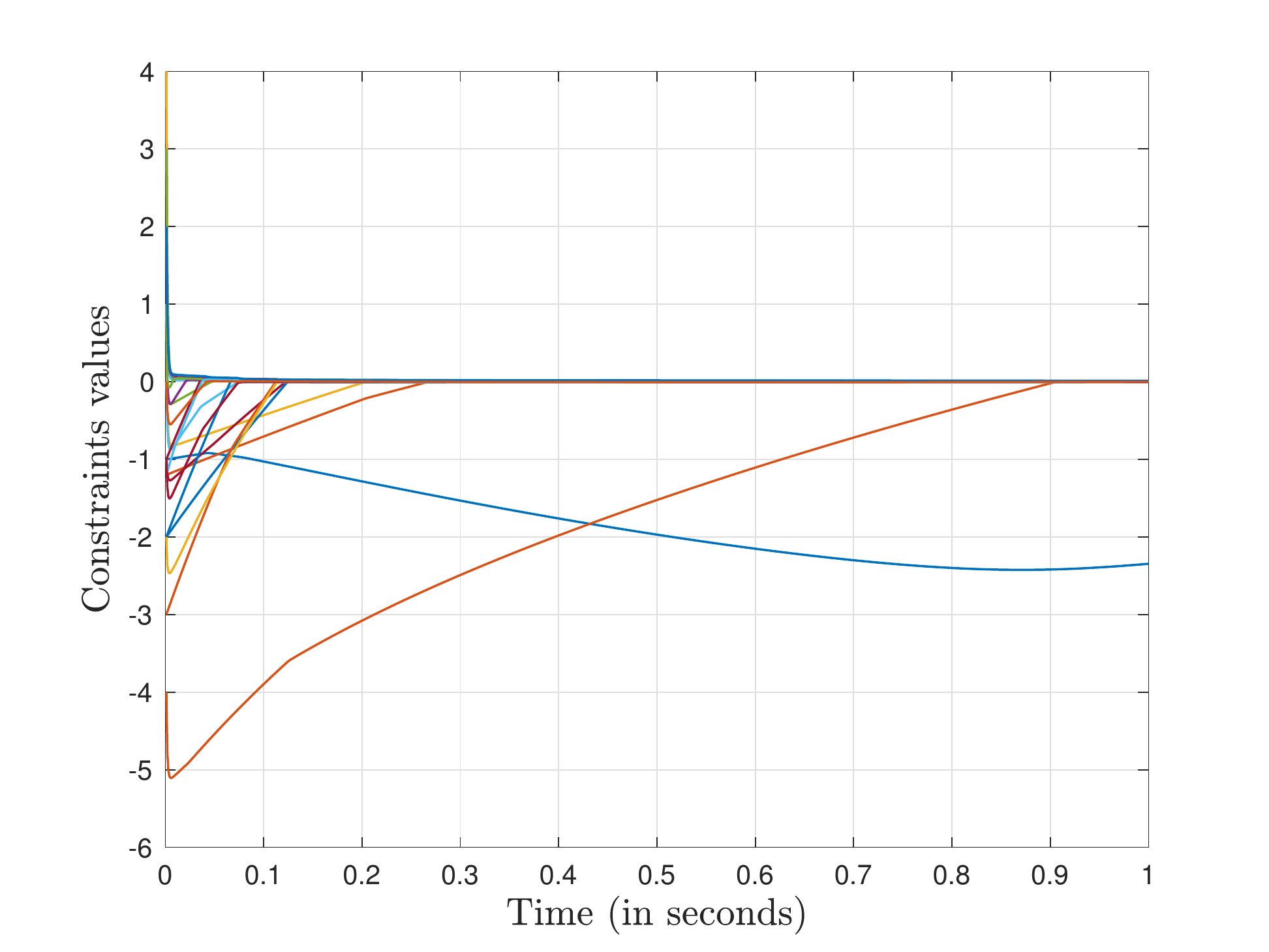}}
  \subfigure[]{\includegraphics[scale=0.45]{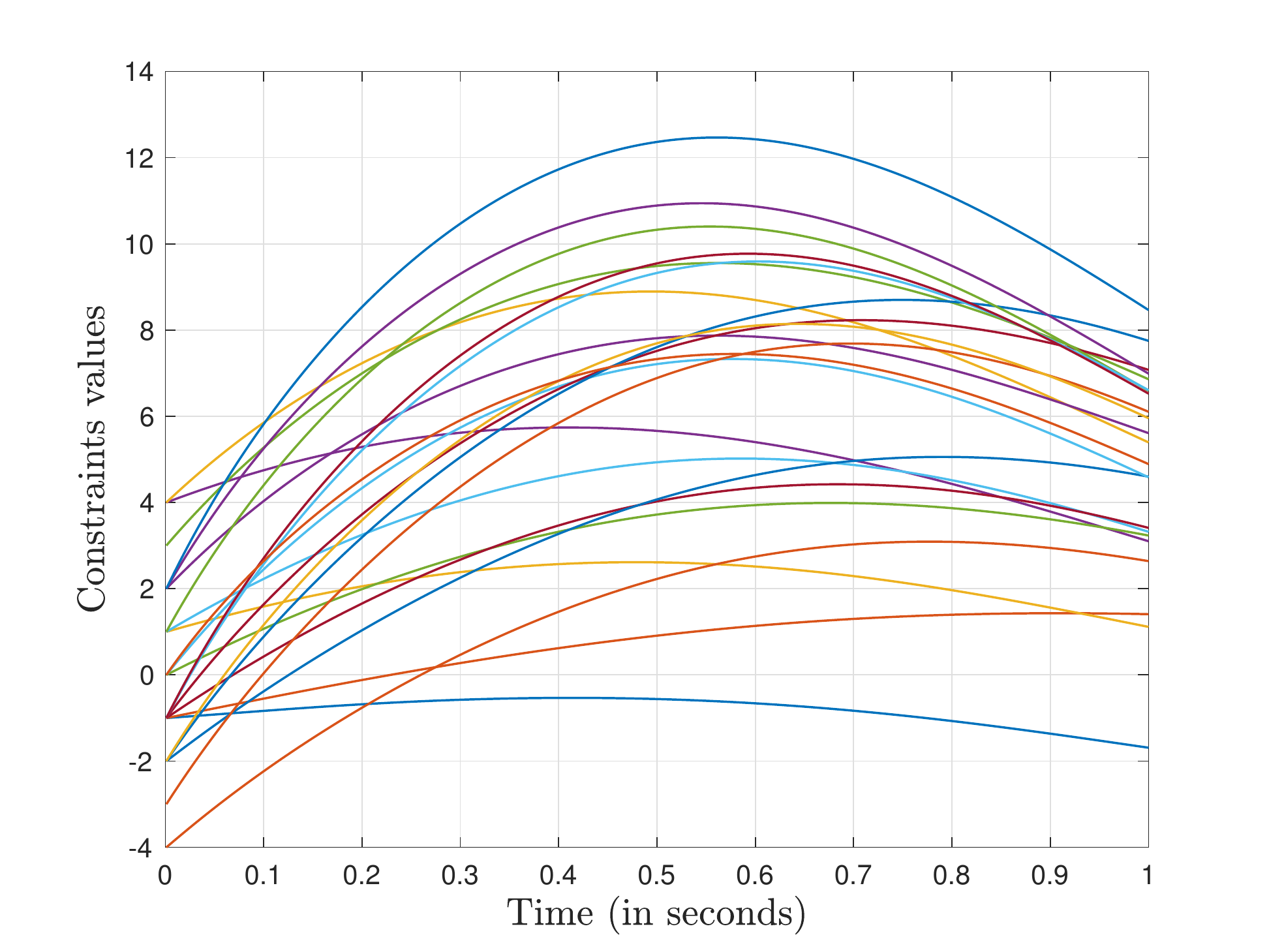}}
  \caption{Evolution of the constraints under (a) the proposed
    distributed dynamics with $\tau=1$ and (b) the saddle-point
    dynamics.}\label{fig:constraints}
\end{figure}

Here, we illustrate the effectiveness of the proposed distributed
dynamics~\eqref{eq:algorithm}. Our optimization problem is inspired
by~\cite{FPK-AKM-DKHT:98}: we consider 50 agents connected in a circle
forming a ring topology and seeking to solve
\begin{equation*}
  \begin{aligned}
    & \max_{x \in \D}
    & & \sum\limits_{i=1}^{50}  f_i(x_i) \\
    &\; \; \text{s.t.}
    & & Ax \leq C .
  \end{aligned}
\end{equation*}
Here, $f_i(x_i)=i \log x_i$ for $i\in \until{50}$.  The sparse matrix
$A \in \real^{23 \times 50}$ is such that each of the 23 constraints
it defines involves a different corresponding agent and its 1-hop
neighbors. We take $\D=\setdef{x \in \real^n}{ 10^{-1} \le \| x
  \|_{\infty} \leq 10}$.  Throughout the simulations, we consider the
exact penalty function~\eqref{eq:penalty} with $\epsilon=10^{-2}$ and
$\gamma=1$.  \revision{Since the dynamics are in continuous time, we
  use a first-order Euler discretization for the MATLAB implementation
  with stepsize $10^{-3}$. We compare the performance of the proposed
  distributed algorithm with values $\tau=1$ and $\tau = 10^{-1}$,
  resp., against the centralized gradient descent~\eqref{eq:idealgd},
  the saddle-point dynamics~\cite{AC-BG-JC:17-sicon} of the
  Lagrangian, and the centralized and the distributed Nesterov's
  accelerated gradient method~\cite{YEN:83} of the penalty function.
  To implement the latter, we use $\tau=1$ and
  replace~\eqref{eq:gdynamics} with Nesterov's acceleration step.  We
  use the same initial condition for all the algorithms.  }
Figure~\ref{fig:comparison} shows the evolution of the objective
function under each algorithm.
One can observe that the proposed distributed algorithm performs much
better than the saddle-point dynamics. As expected, centralized
Nesterov's accelerated gradient method performs the best, followed by
the distributed Nesterov method obtained by applying the acceleration
to our proposed distributed algorithm. The output of the distributed
algorithm for both values of $\tau$ is also close to that of the
centralized gradient descent.  Figure~\ref{fig:constraints} show the
evolution of the value of $Ax-C$ for the proposed distributed
algorithm with $\tau=1$ and the saddle-point dynamics. Even though
Proposition~\ref{prop:feasible} states that, for the centralized
gradient descent counterpart, there is no guarantee of staying inside
the feasible set for general constraints, Figure~\ref{fig:constraints}
shows that the distributed algorithm satisfies the constraints much
better during the evolution than the saddle-point dynamics.

\begin{figure}[htb]
  \centering
  \includegraphics[scale=0.45]{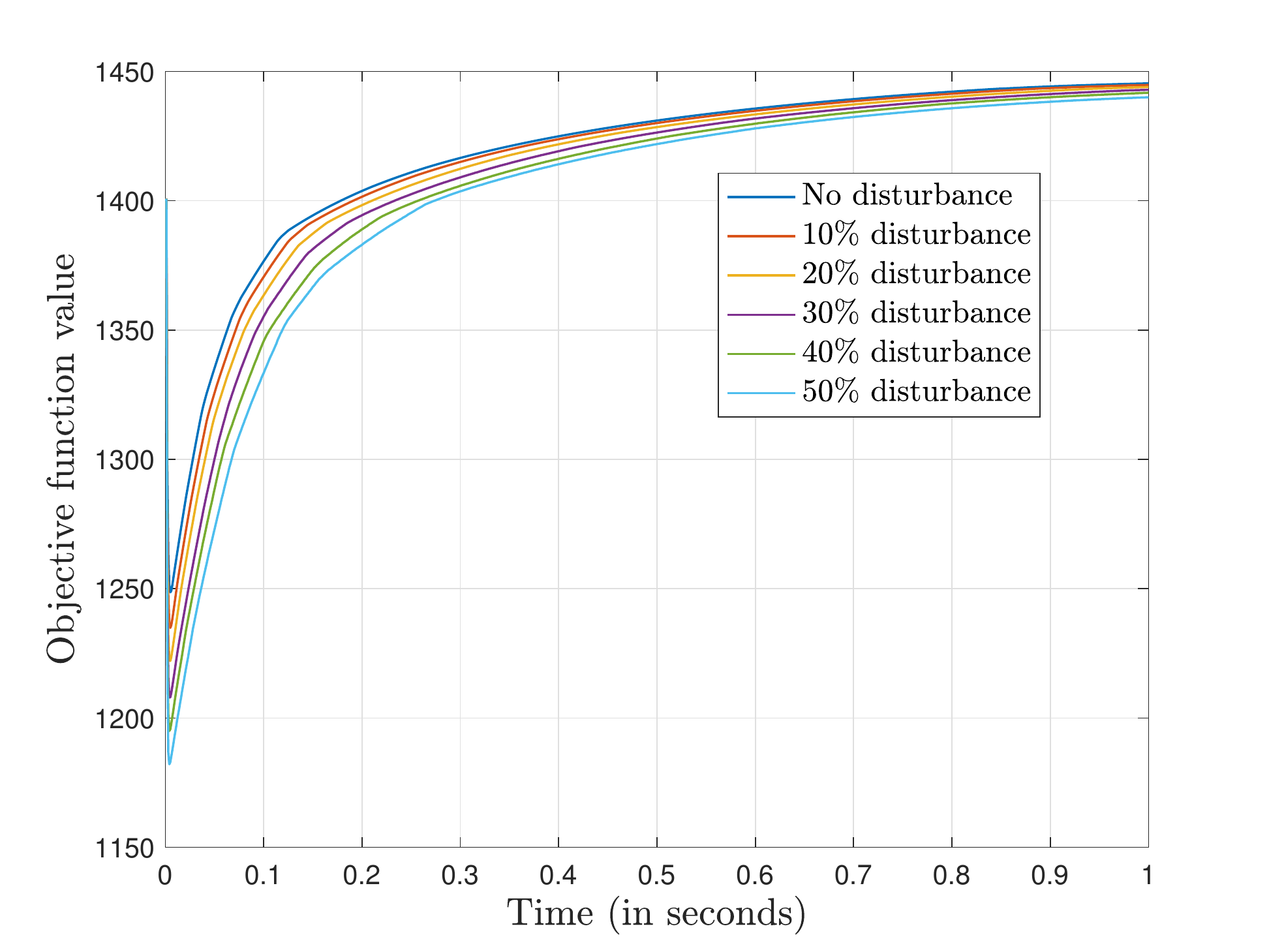}
  \caption{Evolution of the objective function value under the
    proposed distributed dynamics with $\tau=1$ in the presence of disturbances. The
    amount of disturbance in percentage denotes the ratio of
    the norm of the disturbance to the norm of the unperturbed
    dynamics.}\label{fig:disturbed_function}
\end{figure}

In the next simulation we illustrate the robustness of the proposed
dynamics. For this, we add a disturbance to the
dynamics~\eqref{eq:algorithm} using random vectors at each iteration
as follows. For~\eqref{eq:estimation}, we add $d=\beta \| u (x,
\Upsilon) \| \times (\text{unit-norm random vector})$, where we use
the MATLAB function \verb|rand| to generate random numbers between 0
and 1. Similarly, for~\eqref{eq:gdynamics}, we add $d=\beta \| w (x
,\mathcal{P}_{\chi}\Upsilon) \| \times (\text{unit-norm random
  vector})$.
For the scaling constant $\beta$, which also equals the ratio of the
norm of the total disturbance to the norm of the unperturbed dynamics,
we use gradually increasing values between 0.1 to 0.5. 
For each value of $\beta$, we plot the evolution of the
objective function with $\tau=1$ in Figure~\ref{fig:disturbed_function}.  The plot
shows the graceful degradation of the performance as the ratio of the
norm of disturbance to the norm of unperturbed dynamics increases,
demonstrating the effectiveness of the proposed dynamics against
disturbances.

\section{Conclusions}\label{sec:conc}
We have considered the problem of distributed optimization of a
separable function under locally coupled constraints by a group of
agents. Our approach relies on the reformulation of the optimization
problem via a continuously differentiable exact penalty function. To
enable the distributed computation of the gradient of this function,
we have developed a distributed algorithm, of independent interest, to
solve linear algebraic equations defined by separable data. This
algorithm has exponential rate of convergence, is input-to-state
stable, and does not require the individual agent matrices to be
invertible. Building on this, we have introduced dynamics to
asymptotically compute the gradient of the penalty function in a
distributed fashion.  Our algorithmic solution for optimization
consists of implementing gradient descent and Nesterov's accelerated
method with the running estimates provided by this dynamics.  We have
shown the effectiveness of the proposed algorithm in simulation and
compared its performance against a variety of other methods.  Future
work will explore \revision{the design of distributed algorithms for
finding the least-square solutions of
  linear equations defined by separable problem data which only rely
  on 1-hop communication,} distributed ways to determine the timescale
of the estimation dynamics necessary to guarantee convergence, the
characterization of the rate of convergence of the accelerated
implementation, the study of constraint satisfaction along the
executions, and the extension of our approach to problems involving
global, non-sparse constraints.

\appendix

\textit{Proof of Proposition~\ref{prop:feasible}.}  To prove the
result, we examine the Lie derivative of the constraint functions
along the dynamics. We consider the different cases below:

\emph{(Equality constraints):} Given the constraint function $h$,
consider the Lie derivative over the set $\mathcal{F}$,
\begin{align*}
  L_{-\nabla f^\epsilon} h(x)
    &= - \nabla h(x)' \big( \nabla f(x) + \nabla h(x) \mu(x)
  \\
  & \quad + \nabla \mu(x) h(x)+ \frac{2}{\epsilon} \nabla h(x) h(x)
  \big)
  \\
  &= -\nabla h(x)' (\nabla f(x)+ \nabla h(x) \mu(x))
\end{align*}
where we have used the fact that $h(x)=0$ for $x \in \mathcal{F}$.
Substituting the value of $\mu(x)$ from~\eqref{eq:lambda}, 
\begin{align*}
  L_{-\nabla f^\epsilon} h(x) &= -\nabla h(x)' ( \nabla f(x)
  \\
  & \quad - \nabla h(x)N(x)^{-1}\nabla h(x)' \nabla f(x)) = 0.
\end{align*}
This means that the constraint function remains constant along the
gradient dynamics over~$\mathcal{F}$. Hence, $x(t) \in \mathcal{F}$
for all $t \geq t_0$ regardless of the value of $\epsilon$.

\emph{(Scalar inequality constraint):} With only one inequality
constraint defined by a scalar-valued function~$g$, we have $x \in
\mathcal{F}$ iff $g(x) \leq 0$.  To determine the invariance of
the feasibility set, we only need to look at points where $g(x) =
0$.  In this case, the Lie derivative is
\begin{align*}
  L_{-\nabla f^\epsilon} g (x) = -\nabla g(x)' \big(\nabla
  \lambda(x) + \frac{2}{\epsilon} \nabla g(x) \big) y^{\epsilon^2}
  (x) ,
\end{align*}
where we have already used the fact that $g(x)=0$ and the definition
of $\lambda(x)$ from~\eqref{eq:lambda}. Due to LICQ assumption,
$\nabla g(x)' \nabla g(x) > 0$, and $y^{\epsilon^2} (x) \geq 0$. Since
$\nabla \lambda$ is continuous, it is bounded over the compact set
$\mathcal{D}$.  Hence, there exists $\bar{\epsilon}$ such that for all
$\epsilon \in (0, \bar{\epsilon}]$, $ L_{-\nabla f^\epsilon} g (x)
\leq 0$ for all $x$ such that $g(x)=0$.  This means that $x(t) \! \in
\!  \mathcal{F}$ for all $t \geq t_0$.

\emph{(General constraints):} Here we provide a counterexample for the
case with multiple inequality constraints (a similar one can be
constructed for the case of both equality and inequality
constraints). Consider now a vector-valued function $g$.  The
expression of the Lie derivative evaluated at $x$ such that $g(x) = 0$
is
\begin{align*}
  L_{-\nabla f^\epsilon} g = -\nabla g(x)' \Big(\nabla \lambda(x) + \frac{2}{\epsilon}
  \nabla g(x) \Big) Y^\epsilon(x) y^\epsilon (x) .
\end{align*}
The LICQ assumption implies that $\nabla g(x)' \nabla g(x)$ is
positive definite. However, in general, this is not sufficient to
ensure that the trajectory of the gradient dynamics starting from $x$
will remain in $\mathcal{F}$. To see this, consider the following
example.
\begin{equation*}
  \begin{aligned}
    & \min_{x}
    & & (x_1-1)^2 + (x_2+1)^2 \\
    &\; \; \text{s.t.}
    & & x_1 - 6x_2 \leq 0 \\
    &&& -x_1+x_2 \leq 0
  \end{aligned}
\end{equation*}
Take $x=(0;0)$, where $g(x)=0$. After some calculations, it can be
verified that $\lambda(x)=(0;-2)$ and $Y^\epsilon(x) y^\epsilon
(x)=(0;\epsilon)$. As a result, $\nabla g(x)' \nabla g(x)
Y^\epsilon(x) y^\epsilon (x)=(-7\epsilon;2\epsilon)$ and $L_{-\nabla
  f^\epsilon} g=(14;2\epsilon-4)$. The first component of $L_{-\nabla
  f^\epsilon} g$ is independent of $\epsilon$. This means that no
matter what value of $\epsilon$ we choose, $L_{-\nabla f^\epsilon} g
\nleq 0$ when $g(x)=0$.  Hence, the feasible set is not
invariant. $\hfill \blacksquare$
  
  \begin{IEEEbiography}[{\includegraphics[width=1in,height=1.25in,clip,keepaspectratio]{Figures/photo-priyank}}]{Priyank
      Srivastava}
    received the B.Tech degree in electrical engineering from National
    Institute of Technology, Kurukshetra, India in 2012, and the
    M.Tech degree in control \& automation from Indian Institute of
    Technology Delhi, India in 2016. He is currently pursuing Ph.D. in
    mechanical engineering at the University of California San Diego,
    USA.  His current research interests include dynamical systems,
    distributed and fast optimization, and coordination of distributed
    energy resources to enable their participation in energy markets.
\end{IEEEbiography}

\vspace*{-75ex}

\begin{IEEEbiography}[{\includegraphics[width=1in,height=1.25in,clip,keepaspectratio]{Figures/photo-JC}}]{Jorge
    Cort\'{e}s}
  (M'02, SM'06, F'14) received the Licenciatura degree in mathematics
  from Universidad de Zaragoza, Zaragoza, Spain, in 1997, and the
  Ph.D. degree in engineering mathematics from Universidad Carlos III
  de Madrid, Madrid, Spain, in 2001. He held postdoctoral positions
  with the University of Twente, Twente, The Netherlands, and the
  University of Illinois at Urbana-Champaign, Urbana, IL, USA. He was
  an Assistant Professor with the Department of Applied Mathematics
  and Statistics, University of California, Santa Cruz, CA, USA, from
  2004 to 2007. He is currently a Professor in the Department of
  Mechanical and Aerospace Engineering, University of California, San
  Diego, CA, USA. He is the author of Geometric, Control and Numerical
  Aspects of Nonholonomic Systems (Springer-Verlag, 2002) and
  co-author (together with F. Bullo and S.  Mart{\'\i}nez) of
  Distributed Control of Robotic Networks (Princeton University Press,
  2009).  He is a Fellow of IEEE and SIAM. At the IEEE Control Systems
  Society, he has been a Distinguished Lecturer (2010-2014), and is
  currently its Director of Operations and an elected member
  (2018-2020) of its Board of Governors.  His current research
  interests include distributed control and optimization, network
  science, resource-aware control, nonsmooth analysis, reasoning and
  decision making under uncertainty, network neuroscience, and
  multi-agent coordination in robotic, power, and transportation
  networks.
\end{IEEEbiography}

\end{document}